\documentclass[center 15pt]{article}
\usepackage{mathbbold}
\usepackage{mathrsfs}
\usepackage{fancyref}

\usepackage{algorithm}
\usepackage{url}
\usepackage{makecell}
\usepackage{multirow}
\usepackage{algorithmic}

\usepackage{amsmath}
\usepackage{amsfonts}
\usepackage{amssymb}
\usepackage{color}
\usepackage{array}
\usepackage{graphicx}
\usepackage{amsmath}
\usepackage{amssymb}
\usepackage{fancyhdr}
\usepackage{mathrsfs}
\usepackage{mathbbold}
\usepackage{amsfonts}
\usepackage{wasysym}
\usepackage{amscd}
\usepackage{bbm}
\usepackage{color}
\usepackage{array}
\usepackage{graphicx}
\usepackage{graphics}%
\usepackage{hyperref}
\allowdisplaybreaks
\newtheorem{Theorem}{Theorem}[section]

\newtheorem{Remark}[Theorem]{Remark}

\date{}
\title{Superiorization of EM Algorithm and Its Application in Single-Photon Emission Computed Tomography(SPECT)}
\author{Shousheng Luo \\
School of Mathematical Sciences, Peking University,
Beijing 100871, China
\and Tie Zhou\\ School of Mathematical Sciences, Peking University,
Beijing 100871, China}

\begin{document}
\maketitle
\large
\begin{abstract}
In this paper,
we presented an efficient algorithm to implement the regularization reconstruction of SPECT.
Image reconstruction with priori assumptions
is usually modeled as a constrained optimization problem.
However, there is no efficient algorithm to solve it due to the large scale of the problem.
In this paper, we used the
superiorization of the expectation maximization (EM) iteration to
implement the regularization reconstruction of SPECT.
We first investigated the convergent conditions of the EM iteration in the presence of perturbations.
Secondly, we designed the superiorized EM algorithm based on the convergent conditions, and then
proposed a modified version of it.
Furthermore, we gave two methods to generate desired perturbations for two special objective functions.
Numerical experiments for SPECT reconstruction were conducted to validate the
performance of the proposed algorithms.
The experiments show that the superiorized EM algorithms are more stable and robust for
noised projection data and initial image than the classic EM algorithm,
and outperform the classic EM algorithm in
terms of mean square error and visual quality of the reconstructed images.
\end{abstract}
\section{Introduction}
Single-photon emission computed tomography (SPECT), which can visualize the physiological information of various organs with the help of radiopharmaceuticals \cite{ET2004,Herman1980,Natterer2001}, a biochemical molecule labeled with radioactivity.
The gamma-rays emitted by the injected radioactive material are recorded by a gamma camera rotated around the patient.
The goal of SPECT is to reconstruct the radionuclide distribution from the measurements numerically.
 The variety of existing SPECT reconstruction algorithms can be split into a family of analytical
methods \cite{Metz1980,Natterer_a2001,Novikov2002,Noo2001,Noo2007}  and a wide class of iterative techniques \cite{Shepp1982,Osem1994,Hsiao2010}.
\\
\indent
From the analytic point of view, the SPECT reconstruction problem \cite{Natterer2001,Metz1980,Noo2001} is
to invert the attenuated Radon transform(aRt) of $f$(distribution of radiopharmaceutical)
\begin{eqnarray}
R_af(s,\varphi)=\int_Rf(s\theta+t\theta^\perp)e^{-\int_t^\infty \mu(s\theta+\tau\theta^\perp){d\tau}}dt,\label{eq1}
\end{eqnarray}
where $\mu$ is a known function, referred to as the attenuation map of gamma-rays,
$\theta = (\cos\varphi, \sin\varphi)$ and $\theta^\perp = (-\sin\varphi, \cos\varphi)$.
In practice, $f$ and $\mu $ are two functions with compact support $\Omega$. Therefore, the
integrand in (\ref{eq1}) is zero outside a bounded interval, and
this integral is written over $(-\infty, +\infty)$ for convenience.
\\
\indent
For the iterative methods,
the SPECT reconstruction problem is to solve the following linear system \cite{Shepp1982},
\begin{eqnarray}
Ax=b,\label{sys_general1}
\end{eqnarray}
where the elements of the observed data $b = (b_1,b_2\cdots, b_M)^t\in\mathbb{R}^M$, the
unknown image $x = (x_1,x_2\cdots, x_N)^t\in \mathbb{R}^N$
and the system matrix $A = (a_{ij})\in\mathbb{R}^{M\times N}$ are all nonnegative.
Here and in the following, the superscript $^t$ denotes the transpose of vector.
The aim is to reconstruct the unknown $x$ as an image from the projection data $b$ via stable
algorithms. A solution is not feasible with conventional methods
directly because of the noisy projection data $b$, the ill-posedness and  large scale of the problem.
For the SPECT reconstruction, the system matrix $A$  not only can model the attenuation of the gamma-rays,
but also can fuse some realistic factors, such as photon scattering and camera blurring.
\\
\indent
The expectation maximization (EM) algorithm \cite{Shepp1982,Shepp1985} and the algebraic reconstruction technique (ART)\cite{Herman1970,Censor_2001} are two
 widely used technologies in imaging sciences, due to their simplicity, efficiency and performance.
In practice, the EM algorithm is more appropriate for emission tomography including SPECT.
Firstly, the EM algorithm maintains the nonnegative
constraint in the iteration procedure.
Secondly, the EM algorithm is relatively robust against data inconsistencies introduced
by Poisson noise, because it
 seeks to minimize the Kullback-Liebler(K-L) distance between the measured
 data $b$ and the projection of the estimated image $Ax$, which is equivalent to
 maximizing the likelihood of Poisson distribution.
\\
\indent
 For SPECT reconstruction, it is one of the main issues to estimate the radionuclide distribution from low-counts projection data.
This issue occurs quite frequently because of practical constraints, such as imaging hardware and scanning geometry. Furthermore, in order to reduce acquisition time, radiation dose and imaging cost efficiently,
we should decrease  the counts of the projection data. However, this would cause the strong deterioration of the observed data and the under-determinacy ($m\ll n$) of the linear system (\ref{sys_general1}).
In these situations, the reconstructed images by the EM algorithm are usually dominated by various distortions, because the EM algorithm accepts any solution which minimizes the K-L distance and depends on the initial point. \\
\indent
The qualities of the reconstructed images can be improved by regularization reconstruction methods.
Regularization methods are often used techniques to improve the quality of reconstructed images. There are two regularization reconstruction models: unconstraint optimization\cite{Pan2006,pan2008,pan2009,Defrise2011,Pan2012,Shen2012}
\begin{eqnarray}
\min_x\{\phi(x)+{1\over\lambda}{D(b,Ax)}\}\label{eq_nconstr}
\end{eqnarray}
and constraint optimization \cite{Herman2007,Herman2008,Herman2011}
\begin{eqnarray}
\min\phi(x)\text{\quad s.t.\quad} x\in\{x|D(b,Ax)<\epsilon\}\label{eq_con_opt}
\end{eqnarray}
where $\phi$ is a convex function, representing the prior knowledge. $D(\cdot,\cdot)$ denotes a distance function, such as $l^2$ and K-L distances. The parameters $\lambda,\epsilon$ are nonnegative and related to the noise level. The lower the noise level is, the smaller the parameters are.
\\
\indent
For the first optimization problem (\ref{eq_nconstr}), there are different algorithms \cite{pan2009,Defrise2011,Shen2012,Chambolle2011}.
In this paper, we focus on the second optimization problem (\ref{eq_con_opt}).
To our knowledge, there is no optimal algorithm to solve the constraint optimization
problem (\ref{eq_con_opt}) efficiently due to the large scale. In this paper, we resort to an
emerging approach called superiorization \cite{Censor_2010} to implement the regularization reconstruction of SPECT.
\\
\indent
The superiorization of iterative methods, which was first proposed by the authors of \cite{Herman2007},
is a relaxation technology for the constrained optimization problem.
The superiorized algorithm lies between the feasibility-seeking algorithms, which seek a feasible point in the constrained set, and the optimal algorithms,
which seek the minimum point of objective function in the constrained set.
The aim of superiorization algorithm is to
look for a superior instead of the optimal point of the objective function or just a feasible
point in the constrained set.
The basic idea of superiorization is to do the feasibility-seeking algorithms with
perturbations about the objective function.
\\
\indent
The superiorization of
ART algorithms has been studied and applied to the regularization reconstruction of computed tomography(CT)
\cite{Herman2007,Herman2008,Censor_2010,Censor2009,Herman2012}.
The authors of \cite{Herman2007,Censor2009} first investigated the convergence of the two variants of ART under summable perturbations for consistent case. For inconsistent case,
the authors of \cite{Herman2012} proved the convergence of symmetric version of ART in the presence of summable
perturbations.
The superiorization of the EM algorithm was firstly proposed in \cite{Jiang_sup}, and applied to bioluminescence tomography. However, to our knowledge, it is still an open problem about the convergence of the superiorized  EM algorithm.
\\
\indent In this paper, we first discussed the convergence of the EM algorithm in the presence of perturbation in section \ref{th_al}.
A so-called bounded perturbation resilient (BPR) property of ART is vital in the proof of convergence for the
perturbed version of ART.
However, we cannot prove the BPR property of the EM iteration so far,
because of the nonlinearity of the EM operator.
Therefore, we investigated the convergence of the perturbed EM algorithm under the following assumptions.
Firstly, the perturbations should maintain the positivity of iterations. Secondly, the perturbations should
go to zero with the increase of iterates.
Lastly, the perturbed EM iteration should gradually decrease
the K-L distance between the observed data and the projection of estimation by each iteration. \\
\indent
Based on the convergent conditions, we presented the superiorized EM algorithm and its modified version
in section \ref{sup_EM}.
Furthermore, practicable techniques were given to produce ideal perturbations for two special
objective functions, total variation(TV) \cite{ROF1992} and $l^1$-norm \cite{Candes2005},
widely used in imaging sciences.
While the proposed algorithms are applicable to diverse inverse problems,
in this paper we restricted ourselves to demonstrate its usefulness to the SPECT reconstruction.
Numerical results for SPECT reconstruction were given in section \ref{Num_Res}. As our expectation,
the superiorization algorithms output superior image comparing with the classic EM algorithm in terms of MSE and
visual quality. Some conclusions and discussions were given in section \ref{discussion}.
\section{Perturbations Resilience of EM Iteration}\label{th_al}
For the sake of reference, we first introduce some notations and assumptions. In this paper, an image $x$ is
described as a vector of length $N$ with individual elements $x_j$, $j = 1,2,\cdots, N$. When it is necessary to refer to
pixels in the context of a two-dimensional (2D) image we use the double subscript form $x_{p,q}$ , where
\begin{eqnarray}
j = (q-1)W + p, p = 1, 2,\cdots, H, q = 1, 2, \cdots, W,
\end{eqnarray}
and integers $W$ and $H$ are, respectively, the width and height of the 2D image array, which has a total number
of pixels $N = W \times H$. Denote by $\mathbb{R}^N_+$
the region $\{x\geq0|x\in \mathbb{R}^N\}$. For the system matrix $A$,
we set $H_j=\sum_{i=1}^Ma_{ij}$ and $d_i(x)=\sum_{j=1}^Na_{ij}x_j$ for convenience. Furthermore, we introduce
\begin{eqnarray}
g_{ij}(x)&=&\frac{a_{ij}b_i}{d_i(x)},
\end{eqnarray}
and
\begin{eqnarray}
f_j(x)&=&{1\over H_j}\sum_{i=1}^M\frac{a_{ij}b_i}{d_i(x)}={1\over H_j}\sum_{i=1}^Mg_{ij}(x).
\end{eqnarray}
 Let $P$ denote the EM operator, and then the EM iteration $x^{k+1}=P(x^k)=(P_1(x^k),\cdots,P_N(x^k))^t$ for the problem (\ref{sys_general1}) is defined as
\begin{eqnarray}
x^{k+1}_j&=&P_j(x^k)={x^k_j\over H_j}\cdot\sum\limits_{i=1}^M{b_ia_{ij}\over{\sum_{t=1}^Na_{it}x_t^k}}=x_j^k\cdot f_j(x^k),\label{EqEM1}
\end{eqnarray}
with an initial image $x^0>0$.
Similarly, the perturbed version of the EM iteration is defined as
\begin{eqnarray}
x^{k+1}_j= P_j(x^k + \beta_kv^k) ={y^k_j\over H_j}\cdot\sum_{i=1}^M\frac{b_ia_{ij}}{\sum_{t=1}^Na_{it}y^k_t}=y^k_j\cdot f_j(y^k),\label{eq_perturb1}
\end{eqnarray}
where $y^k = x^k+\beta_kv^k$. Here, the vector $v^k\in \mathbb{R}^N$ and number $\beta_k\geq0$
represent the direction and length of the perturbation of the $k$th iteration, respectively. Hereafter, we call the EM iteration without perturbation as the classic EM iteration to distinguish from each other.
\\
\indent It has been established that the sequence $\{x^n\}$ generated by the classic EM iteration converges to
a minimizer of the K-L distance $I_A^b(x)$ between $b$ and $Ax$ on $R_+^N$, where $I_A^b(x)$ is defined as
\begin{eqnarray}
I_A^b(x)&=&I(b,Ax)=\sum\limits_{i=1}^Mb_i\ln\frac{b_i}{d_i(x)}-\sum\limits_{i=1}^M(b_i-d_i(x)),\label{Eqdd}
\end{eqnarray}
where $I(\cdot,\cdot)$ denote the K-L distance function of any two nonnegative vectors.
 An important inequality
 used in this paper is
\begin{eqnarray}
\ln t\geq 1-{1\over t},\label{ineq_imp}
\end{eqnarray}
for $t>0$, and the inequality holds with equality if and only if $t=1$.
By using the inequality (\ref{ineq_imp}), we have that
$I(x,y)\geq0$ for all vectors $x,y\geq0$, and the inequality holds with equality if and only if $x=y$.
\\
\indent
 Obviously, if $\beta_k=0$, the perturbed EM iteration is the same as  the classic EM iteration,
 in which we are not interested. Therefore, in the following we assume $\beta_k>0$.
A natural question is that under what assumptions the sequence $\{x^n\}$ generated by the perturbed EM
iteration (\ref{eq_perturb1}) converges to a minimizer of (\ref{Eqdd}) as well.
Before discussing the convergent conditions of the perturbed EM iteration,
we first summarize some propositions of the classic EM iteration \cite{Shepp1982,Shepp1985,Deblurring1992}.
 %and {\red K-L distance }.
Without loss of generality, we assume that all the elements $b_i>0$ for all $i$, and
$H_j=\sum_{i=1}^Ma_{ij}=1$ for all $j$ for simplicity.
\begin{Theorem}\label{PropEM}
For any initial image $x^0>0$, denote by $x^k$ the estimate of the classic EM algorithm after $k$ iterations.
Then the following propositions hold:
\begin{enumerate}
\item $x^k>0$ and $\sum_{j=1}^Nx_j^k=\sum_{i=1}^Mb_i$, for all $k>0$.
\item\label{item2}
$I_A^b(x^{k+1})\leq I_A^b(x^k)$, and $I(x^{k+1},x^{k})\leq I_A^b(x^{k})-I_A^b(x^{k+1})$.
\item $\{x^k\}$ converges to a minimizer $x^\ast$ of $I_A^b(x)$ on $\mathbb{R}_+^N$, and $x^\ast$ is a fixed
point of the EM operator.
\item \label{item5}$I(x^\ast,x^{k+1})\leq I(x^\ast,x^k)$.
\end{enumerate}
The inequalities in second and fourth items above hold with equalities if and only if $x^k$ is a fixed point of the EM operator.
\end{Theorem}
%%%%%%%%%%%%%%%%%%%%%%%%%%%%%%%%
%%%%
\begin{Remark}
Obviously, $x$ is a fixed point of the EM operator if and only if $f_j(x)=1$ for $x_j\neq 0$.
\end{Remark}
From the propositions above, we have that the sequence
$\{I_A^b(x^k)\}$ monotonically converges to the minimum of $I_A^b(\cdot)$ on $\mathbb{R}_+^N$,
and $\{x^k\}$ approximates to the minimizer $x^\ast$ gradually.
Next, we investigate the convergence of the perturbed EM iteration, and prove the similar propositions as far as we can.
\begin{Theorem}\label{ThEM}
Given any initial image $x^0$, denote by $\{x^k\}_{k\in\mathbb{N}}$ the sequence generated by the perturbed EM iteration (\ref{eq_perturb1}).
\begin{eqnarray}
x^{k+1}_j&=&({x_j^k+\beta_k{v^k_j}})\cdot f_j(x^k+\beta_kv^k)
={{y^k_j}}\cdot f_j(y^k),\label{eq_pEM}
\end{eqnarray}
where $y_j^k=x^k_j+\beta_kv_j^k$, and $\{v^k\}$ is a bounded sequence.
Suppose that
\begin{enumerate}
\item Positivity: $y_j^k>0$.
\item  Vanishing: $\beta_k \longrightarrow 0$.
\item Decreasing
\begin{eqnarray}
\beta_k\max_{j\in{S_k^-}}\{-\frac{v_j^k}{y_j^k}\}B_{k}^--\beta_k\min_{j\in{S_k^+}}\{\frac{v_j^k}{y_j^k}\}B_{k}^+
+\beta_k\sum_{j=1}^Nv_j^k&<&I_A^b(y^k)-I_A^b(x^{k+1}),\label{ieq_1}
\end{eqnarray}
where
$
S^-_k=\{j|v_j^k<0\},S^+_k=\{j|v_j^k>0\}
$
, and
$B_k^-=\max\{\rho+\sum_{j\in{S_k^-}}x_j^{k+1},\sum_{j\in{S_k^-}} \hat{x}_j\}$,
$B_k^+=\min\{\sum_{j\in{S_k^+}}x_j^{k+1},\sum_{j\in{S_k^+}}\hat{x}_j\}$.
Here $\rho$ is a sufficiently small positive number.
\end{enumerate}
Then, the following propositions hold:
\begin{enumerate}
\item $I_A^b(x^{k+1})\leq I_A^b(x^{k})$.
\item $\{x^k\}$ has a convergent subsequence $\{x^{m_k}\}$, and the limit $\hat{x}$  is a fixed point of the EM operator.
\item  $I(\hat{x},x^{k+1})\leq I(\hat{x},x^k)$, and $x^k\longrightarrow \hat{x}$.
\item  $\hat{x}$ is a minimizer of $I_A^b(x)$.
\end{enumerate}
Furthermore, the inequalities above hold with equalities iff $x^k$ is a fixed point of the EM operator and $\beta_kv^k=0$.
\end{Theorem}
\indent Before presenting the proof, we explain the necessities of the conditions of
theorem \ref{ThEM}. The positive condition of $y_j^k$ is necessary for the nonnegative constraints.
The second condition is required by the convergence of the sequence $\{x^k\}$.
Intuitively, the last condition is used to guarantee that $I_A^b(x^{k+1})\leq I_A^b(x^{k})$,
which implies the convergence of $I_A^b(x^k)$.
\\
\indent An important concern is about the existence of the perturbations satisfying the conditions of theorem \ref{ThEM}, especially the third condition.
Because the sequence $\{\hat{x}^k\}$ generated by the classic EM iteration ($\beta_k=0$) satisfies the conditions in theorem \ref{ThEM}, the sequence $\{x^k\}$ generated by the perturbed version also satisfies them when $\beta_k$s are small enough for each iteration due to the continuities of the EM operator (\ref{eq_perturb1}) and K-L distance (\ref{Eqdd}). Therefore, there exist
perturbations which satisfy the assumptions of theorem \ref{ThEM}.\\
{\bf Proof : } We prove this theorem step by step, which follows the proving procedure in
\cite{Deblurring1992}.
\begin{enumerate}
%%%%%%%%%%%%%%%%%%%%%%
\item {\bf Proof of proposition 1:}
By the definition of $I_A^b(\cdot)$ and $d_i(\cdot)$, we have
\begin{eqnarray}
& &I_A^b(x^k)-I_A^b(x^{k+1})\notag\\
&=&\sum_{i=1}^Mb_i\ln\frac{b_i}{d_i(x^k)}-\sum_{i=1}^Mb_i\ln\frac{b_i}{d_i(x^{k+1})}\notag\\
&=&\sum_{i=1}^Mb_i\ln\frac{d_i(x^{k+1})}{d_i(y^k)-\beta_kd_i(v^k)}\notag\\
&=&\sum_{i=1}^Mb_i\ln\frac{d_i(x^{k+1})}{d_i(y^k)}+\sum_{i=1}^Mb_i\ln\frac{1}{1-\beta_kd_i(v^k)/d_i(y^k)}\notag\\
&\geq&I_A^b(y^k)-I_A^b(x^{k+1})-\beta_k\sum_{j=1}^Nv_j^k+\sum_{i=1}^Mb_i\left[1-\left(1-\beta_k\frac{d_i(v^k)}
{d_i(y^k)}\right)\right]\label{ieq_23}\\
&=&I_A^b(y^k)-I_A^b(x^{k+1})-\beta_k\sum_{j=1}^Nv_j^k+\sum_{i=1}^Mb_i\beta_k\frac{d_i(v^k)}{d_i(y^k)}\notag\\
&\geq&\beta_k\max_{j\in{S_k^-}}\{-\frac{v_j^k}{y_j^k}\}B_{k}^--\beta_k\min_{j\in{S_k^+}}\{\frac{v_j^k}{y_j^k}\}B_k^+
+\beta_k\sum_{i=1}^Mb_i\frac{d_i(v^k)}{d_i(y^k)}\notag\\
&=&\beta_k\max_{j\in{S_k^-}}\{-\frac{v_j^k}{y_j^k}\}B_k^--\beta_k\min_{j\in{S_k^+}}\{\frac{v_j^k}{y_j^k}\}B_k^++\beta_k\sum_{j=1}^N\frac{v_j^k}{y_j^k}y_j^k\sum_{i=1}^M\frac{a_{ij}b_i}{d_i(y^k)}\notag\\
&=&\beta_k\max_{j\in{S_k^-}}\{-\frac{v_j^k}{y_j^k}\}B_k^--\beta_k\min_{j\in{S_k^+}}\{\frac{v_j^k}{y_j^k}\}B_k^++\beta_k\sum_{j=1}^N\frac{v_j^k}{y_j^k}x^{k+1}_j\notag\\
&\geq&\beta_k\max_{j\in{S_k^-}}\{-\frac{v_j^k}{y_j^k}\}(B_k^--\sum_{j\in{S_k^-}}x^{k+1}_j)-\beta_k\min_{j\in{S_k^+}}\{\frac{v_j^k}{y_j^k}\}(B_k^+-\sum_{j\in{S_k^+}}x^{k+1}_j)\notag\\
&\geq&\beta_k\max_{j\in{S_k^-}}\{-\frac{v_j^k}{y_j^k}\}\rho-\beta_k\min_{j\in{S_k^+}}\{\frac{v_j^k}{y_j^k}\}(B_k^+-\sum_{j\in{S_k^+}}x^{k+1}_j).\label{min_use}
\end{eqnarray}%%%%%%%%%%%%%%%%%%%%%%%%%
The first and second inequalities hold because of the inequality (\ref{ineq_imp}) and the condition (\ref{ieq_1}), respectively.
Therefore, we have proved that $I_A^b(x^{k+1})\leq I_A^b(x^k)$ by  (\ref{ieq_1}) and (\ref{min_use}).
\\
\indent
Finally, we should prove that  $I_A^b(x^{k+1})=I_A^b(x^k)$ iff $\beta_kv^k=0$ and $x^k$ is a fixed point of EM operator.
The sufficiency is obvious by theorem \ref{PropEM}.
The necessity is also true by the following facts.  If $I_A^b(x^{k+1})=I_A^b(x^k)$, we have that all the inequalities in the derivation above hold with equalities, which implies that $\frac{1}{1-\beta_{k}d_i(v^k)/d_i(y^k)}=1$ for all $i$ and $S^-_k=\emptyset$ by inequalities (\ref{ieq_23}) and (\ref{min_use}), respectively,
i.e. $\beta_kd_i(v^k)=0$ and $v_j^k\geq0$. In fact, we have that $v_j^k=0(\forall~j)$, i.e. $y_j^k=x_j^k$, and  $x^k$ is a fixed point by theorem \ref{PropEM}. Otherwise, assume $v_{j_0}^k>0$ for some $j_0$, there must exist $i_0$ such that $a_{i_0j_0}>0$ and
$d_{i_0}(v^k)=\sum_{j=1}^Na_{i_0j}v^k_j\geq a_{i_0j_0}v_{j_0}>0$ by
the assumption $\sum_{i=1}^Ma_{ij}=1\neq0$, which is contradict to $d_i(v^k)=0(\forall~i)$.
\item {\bf Proof of proposition 2: } Because
$0<\sum_{j=1}^Nx_j^k=\sum_{i=1}^Mb_i$(constant) and $x^k>0$, $\{x^k\}$ has a
convergent subsequence $\{x^{m_k}\}$. Denoting by $\hat{x}$ the limit
of $\{x^{m_k}\}$,
we have $y^{m_k}\longrightarrow \hat{x}$ because $\beta_k\longrightarrow 0$ and $v^k$ is bounded.
%%%%%%%%%%%%%%%%%%%%%%%%%%%%%%%%%%%%%%%%%%%%%%%%%%%55555555555555555555555555555
Next we  prove $\hat{x}$ is a fixed point of the EM operator. To this end,
we define a function for $x\geq0$
\begin{eqnarray}
D(x)=I(P(x),x).
\end{eqnarray}
By the second proposition of theorem \ref{PropEM}, we have
\begin{eqnarray}
D(y^{m_k})=I(P(y^{m_k}),y^{m_k})=I(x^{m_k+1},y^{m_k})\leq I_A^b(y^{m_k})-I_A^b(x^{m_k+1}).
\end{eqnarray}
Due to the proposition 1 of theorem \ref{ThEM}, $\{I_A^b(x^k)\}$ is a convergent sequence. Furthermore, $\{I_A^b(y^k)\}$ converges to the same  limit as $\{I^b_A(x^k)\}$ because of the continuity of $I^b_A(\cdot)$ and $\beta_k\longrightarrow0$.
Therefore, we have $I_A^b(y^{m_k})-I_A^b(x^{m_k+1})\longrightarrow 0$, $m_k\longrightarrow \infty$.
Thus, we have that
\begin{eqnarray}
D(\hat{x})=\lim_{m_k\longrightarrow\infty}D(y^{m_k})=I(P(\hat{x}),\hat{x})=0,
\end{eqnarray}
i.e. $\hat{x}$ is a fixed point of the EM operator, since $I(x,y)=0\Longleftrightarrow x=y$.
Therefore, we have $f_j(\hat{x})=1$ if $\hat{x}_j\neq0$ by the EM iteration formula.
%%%%%%%%%%%%%%%%%%%%%%%%%%%%%%%%%
%%%%%%%%%%%%%%%%%%%%%%%%%%%%
\item {\bf Proof of proposition 3: }
 Before going further,
we first prove a general inequality that $I(\hat{x},x)-I(\hat{x},P(x))\geq I_A^b(x)-I_A^b(\hat{x})$ for
any point $x>0$. Assuming $z=P(x)$, we have
\begin{eqnarray}
& &I(\hat{x},x)-I(\hat{x},z)\notag\\
&=&\sum_{j=1}^N\hat{x}_j\ln\frac{\hat{x}_j}{x_j}-\sum_{j=1}^N(\hat{x}_j-x_j)+
\sum_{j=1}^N\hat{x}_j\ln\frac{\hat{x}_j}{z_j}-\sum_{j=1}^N(\hat{x}_j-z_j)\notag\\
&=&\sum_{j=1}^N\hat{x}_j\ln\frac{z_j}{x_j}-\sum_{j=1}^N(\hat{x}_j-x_j)\notag\\
&=&\sum_{j=1}^N\hat{x}_j\sum_{i=1}^Mg_{ij}(\hat{x})\ln\frac{z_j}{x_j}\frac{g_{ij}(x)}{g_{ij}(x)}\frac{g_{ij}(\hat{x})}{g_{ij}(\hat{x})}-\sum_{j=1}^N(\hat{x}_j-x_j)\notag\\
&=&\sum_{j=1}^N\hat{x}_j\sum_{i=1}^Mg_{ij}(\hat{x})\ln\frac{z_j}{x_j}\frac{g_{ij}(\hat{x})}{g_{ij}(x)}+\sum_{j=1}^N\hat{x}_j\sum_{i=1}^Mg_{ij}(\hat{x})\ln\frac{g_{ij}(x)}{g_{ij}(\hat{x})}-\sum_{j=1}^N(\hat{x}_j-x_j)\notag\\
&\geq&-\sum_{j=1}^N\hat{x}_j\ln\sum_{i=1}^Mg_{ij}(\hat{x})\frac{x_j}{z_j}\frac{g_{ij}(x)}{g_{ij}(\hat{x})}+\sum_{j=1}^N\hat{x}_j\sum_{i=1}^Mg_{ij}(\hat{x})\ln\frac{\frac{a_{ij}b_i}{d_i(x)}}{\frac{a_{ij}b_i}{d_i(\hat{x})}}-\sum_{j=1}^N(\hat{x}_j-x_j)\notag\\
&=&-\sum_{j=1}^N\hat{x}_j\ln\frac{x_jf_j(x)}{z_j}+\sum_{i=1}^M\sum_{j=1}^N\hat{x}_j\frac{a_{ij}b_i}{d_i(\hat{x})}\ln\frac{d_i(\hat{x})}{d_i(x)}-\sum_{j=1}^N(\hat{x}_j-x_j)\notag\\
&=&\sum_{i=1}^Mb_i\ln\frac{d_i(\hat{x})}{d_i(x)}-\sum_{j=1}^N(\hat{x}_j-x_j)\notag\\
&=&\sum_{i=1}^Mb_i\ln\frac{b_i}{d_i(x)}-\sum_{i=1}^M(b_i-d_i(x))+\sum_{i=1}^M(b_i-d_i(x))\notag\\
&&-\sum_{i=1}^Mb_i\ln\frac{b_i}{d_i(\hat{x})}
+\sum_{i=1}^M(b_i-d_i(\hat{x}))-\sum_{i=1}^M(b_i-d_i(\hat{x}))-\sum_{j=1}^N(\hat{x}_j-x_j)\notag\\
&=&I_A^b(x)-I_A^b(\hat{x})+\sum_{i=1}^M(d_i(\hat{x})-d_i(x))-\sum_{j=1}^N(\hat{x}_j-x_j)\notag\\
&=&I_A^b(x)-I_A^b(\hat{x}),\label{eq_31}
\end{eqnarray}
\large
where we used the relationship $\sum_{j=1}^N\hat{x}_j=\sum_{j=1}^Nx^k_j=\sum_{i=1}^Mb_i=\sum_{j=1}^Nz_j$.
The inequality and the second equality are implied by Jensen's inequality and the fact $f_j(\hat{x})=\sum_{i=1}^Mg_{ij}(\hat{x})=1$ for $\hat{x}_j\neq 0$, respectively.
The last equality follows the fact that
\begin{eqnarray}
\sum_{i=1}^Md_i(x)=\sum_{i=1}^M\sum_{j=1}^Na_{ij}{x}_j=\sum_{j=1}^Nx_j,
\end{eqnarray}
 by
 the assumption that $\sum_{i=1}^Ma_{ij}=1$.
Let $x=y^k$, then $z=P(y^k)=x^{k+1}$ and $I(\hat{x},y^k)-I(\hat{x},x^{k+1})\geq I_A^b(y^k)-
I_A^b(\hat{x})$ by (\ref{eq_31}).
Furthermore, we have
\begin{eqnarray}
I(\hat{x},y^k)-I(\hat{x},x^{k+1})
&\geq&I_A^b(y^k)-I_A^b(\hat{x})\notag\\
&=&\left[I_A^b(y^k)-I_A^b(x^{k+1})\right]+I_A^b(x^{k+1})-I_A^b(\hat{x})\notag\\
&\geq&I_A^b(y^k)-I_A^b(x^{k+1}).\label{eq_cond}
\end{eqnarray}
The last inequality holds because $I_A^b(x^{k+1})\geq I_A^b(\hat{x})$.
Now, we prove  $I(\hat{x},x^{k+1})\leq I(\hat{x},x^{k})$.
By the definitions of the K-L distance and the perturbed EM iteration, we have
%TTTTTTTTTTTTTT%%%%%%%%%%%%%%%%%%%%%%%%%%%%%%%%%%%%%%%%%%%%%%%%%%%%%%%%%%%%%%%%%%%%%
\begin{eqnarray}
& & I(\hat{x},x^k)-I(\hat{x},x^{k+1})\notag\\
&=&\sum_{j=1}^N\hat{x}_j\ln\frac{(x_j^k+\beta_kv^k_j)f_j(y^k)}{x^k_j}\notag\\
&=&\sum_{j=1}^N\hat{x}_j\ln{f_j(y^k)}+\sum_{j=1}^N\hat{x}_j\ln(1+\beta_k\frac{v_j^k}{x_j^k})\notag\\
&\geq&\sum_{j=1}^N\hat{x}_j\ln\frac{y^k_jf_j(y^k)}{y_j^k}+\sum_{j=1}^N\hat{x}_j\left(1-\frac{1}{1+\beta_k{v_j^k}/{x_j^k}}\right)\label{ieq_25}\\
&=&I(\hat{x},y^k)-I(\hat{x},x^{k+1})+\sum_{j=1}^N\hat{x}_j\frac{\beta_kv_j^k}{x_j^k+\beta_k{v_j^k}}-\beta_k\sum_{j=1}^Nv_j^k\notag\\
&\geq&I_A^b(y^k)-I_A^b(x^{k+1})-\beta_k\max_{j\in{S_k^-}}\{-\frac{v_j^k}{y_j^k}\}\sum_{j\in{S_k^-}}\hat{x}_j+\beta_k\min_{j\in{S_k^+}}\{\frac{v_j^k}{y_j^k}\}\sum_{j\in{S_k^+}}\hat{x}_j-\beta_k\sum_{j=1}^Nv_j^k\notag\\
&\geq&\beta_k\max_{j\in{S_k^-}}\{-\frac{v_j^k}{y_j^k}\}(B_k^--\sum_{j\in{S_k^-}}\hat{x}_j)+\beta_k\min_{j\in{S_k^+}}\{\frac{v_j^k}{y_j^k}\}(\sum_{j\in{S_k^+}}\hat{x}_j-B_k^+)\notag\\
&\geq&0.\label{ieq_cond1}
\end{eqnarray}
Again we used the inequality (\ref{ineq_imp}) for the first inequality.
 The second and the last inequalities hold by (\ref{eq_cond}) and (\ref{ieq_1}), respectively.
Thus, we have that $I(\hat{x},x^{k+1})\leq I(\hat{x},x^{k})$ by equation (\ref{ieq_cond1}).\\
 Next, we need to prove that
$I(\hat{x},x^{k+1})= I(\hat{x},x^{k})$ iff $\beta_kv^k_j=0$ and $x^k$ is a fixed point of EM operator. The sufficiency is clear.
The necessity is also true by the following facts.  If $I(\hat{x},x^{k+1})= I(\hat{x},x^{k})$, we have that all the inequalities in the derivation above hold with equalities, which implies that $\beta_kv_j^k=0$ and $I_A^b(y^k)=I_A^b(x^{k+1})$ by inequalities (\ref{ieq_25}) and (\ref{ieq_cond1}), respectively.
Therefore, we have that $y^k=x^k$ and $x^k$ is a fixed point of EM operator by theorem \ref{PropEM}.
\\
\indent
Finally, we need to prove that $x^k\longrightarrow \hat{x}$.
Firstly, we have $I(\hat{x},x^k)\searrow0$ since $I(\hat{x},x^{m_k})\longrightarrow0$ and
the sequence $\{I(\hat{x},x^k)\}$ monotonously decreases. Therefore,
we have proved $x^k\longrightarrow \hat{x}$ by
$I(x,y)=0\Longleftrightarrow x=y$, and $y^k\longrightarrow \hat{x}$ as well by $\beta_kv^k\longrightarrow 0$.
\item {\bf Proof of proposition 4: }
In order to prove this proposition, we should prove that $\hat{x}$ satisfies the Kuhn-Tucker (K-T) conditions. For the EM algorithm and the K-L function,
the K-T conditions are equivalent to (1) $f_j(\hat{x})=1$ if $\hat{x}_j\neq0$, and
(2) $0\leq f_j(\hat{x})\leq1$ if $\hat{x}_j=0$ \cite{Shepp1982,Deblurring1992}. Obviously, the first condition holds because
$\hat{x}$ is a fixed point of the EM operator. For the second condition, the nonnegativity of $f_j(\hat{x})$ is satisfied since $x^k_j\longrightarrow\hat{x}_j\geq0$, $b_i\geq0$ and $a_{ij}\geq0$.
Next we need to prove that $f_j(\hat{x})\leq1$ for $\hat{x}_j=0$. Suppose $f_j(\hat{x})\geq1+\epsilon>1$ and  $\hat{x}_j=0$ for some $j$.
There exists a sufficiently large integer $L>0$
such that $f_j(y^k)>1+{\epsilon_1}$ with $\epsilon_1>0$ for all $k\geq L$
since $y^k\longrightarrow \hat{x}$.
If $x_j^k\longrightarrow 0$, there must be infinite iterations such that $x^{k+1}_j<x^k_j$. In the following, we assume that $k>L$.
For each iteration satisfying $x^{k+1}_j=y_j^k\cdot{f_j(y^k)}<x^k_j$, we have that $y_j^k<x_j^{k+1}<x_j^k$ since $f_j(y^k)>1$, and
\begin{eqnarray}
1&>&\frac{(x^k_j+\beta_kv_j^k)f_j(y^k)}{x^k_j}\notag\\
 &=&(1+\beta_k\frac{v_j^k}{x_j^k})f_j(y^k).
\end{eqnarray}
By the assumption that $f_j(y^k)>(1+\epsilon_1)$, we have
\begin{eqnarray}
1+\beta_k\frac{v_j^k}{x_j^k}<\frac{1}{1+\epsilon_1}.\label{eq_uniq}
\end{eqnarray}
Abstracting 1 on both sides of (\ref{eq_uniq}), we have $\beta_k\frac{v_j^k}{x_j^k}<-\frac{\epsilon_1}{1+\epsilon_1}$.
By equation (\ref{min_use}), we have that
\begin{eqnarray}
& &I_A^b(x^k)-I_A^b(x^{k+1})\notag\\
&\geq&\beta_k\max_{j\in{S_k^-}}\{-\frac{v_j^k}{y_j^k}\}(B_k^--\sum_{j\in{S_k^-}}x^{k+1}_j)-\beta_k\min_{j\in{S_k^+}}\{\frac{v_j^k}{y_j^k}\}(B_k^+-\sum_{j\in{S_k^+}}x^{k+1}_j)\notag\\
&\geq& \frac{\epsilon_1}{1+\epsilon_1}(B_k^--\sum_{j\in{S_k^-}}x^{k+1}_j)\notag\\
&\geq&\rho\frac{\epsilon_1}{1+\epsilon_1}.\label{min_used2}
\end{eqnarray}
The second inequality holds since $y_j^k<x_j^k$ and $B_k^+\leq\sum_{j\in{S_k^+}}x^{k+1}_j$. The last inequality follows from the definition of $B_k^-$.
\\
\indent
From equation (\ref{min_used2}), we can see that the decrease of $I_A^b(x^k)$ is larger than a positive number for each iteration satisfying $x_j^{k+1}<x_j^k$. Furthermore, because the number of iterations satisfying $x_j^{k+1}<x_j^k$ is infinite and $I_A^b(x^k)$ is finite,
we have $I_A^b(x^l)\longrightarrow -\infty$, which is contradict to the fact that $I_A^b(x)\geq0$ for all $x\geq0$.
\end{enumerate}
\section{Superiorization of the EM Algorithm}\label{sup_EM}
The SPECT regularization reconstruction can be modeled as a
constrained optimization problem
\begin{eqnarray}
&&\min_{x\in E}\phi(x),\textrm{\quad\quad } E=\{x^\ast|x^\ast=\text{arg}\min_{x\geq0}I_A^b(x)\},\label{prob1}
\end{eqnarray}
where $\phi$ is a convex function,
which assigns each image $x$ a number indicating the "undesirability" of the image in some sense.
The set $E$ is called feasible set, and it is called feasible problem to look for a point in $E$.
\\
\indent
To our knowledge, there is no efficient algorithm to deal with the constrained optimization problem (\ref{prob1}) because of the large scale of it for SPECT reconstruction. On the other hand,
although the feasible problem is also a constrained optimization, we can solve it by the classic EM algorithm \cite{Shepp1982} and its variants
\cite{Osem1994,Hsiao2010} efficiently. Based on the facts above,
we use the superiorization methodology
to implement the SPECT regularization reconstruction.
\\
\indent  For the objective function $\phi$, the superiorized
 EM algorithm is illustrated as algorithm \ref{algm1} based on the conditions of theorem \ref{ThEM}.
 In order to emphasize the objective function $\phi$ for which we are superiorizing, we refer to the superiorized algorithm as the $\phi$-superiorization of the EM iteration.

 \begin{algorithm}[h]
\caption{ Framework of $\phi$-superiorization algorithm} %算法的标题
\label{algm1} %给算法一个标签，这样方便在文中对算法的引用
\begin{algorithmic}
\STATE Initialization: $\beta_0>0$, $x^0>0$, $k=0$, and $0<\gamma<1$.
\STATE repeat: logic=true
\STATE \indent~~while logic
\STATE \indent~~~~ find a decreasing direction $v^k$ of $\phi$ at $x^k$, such that
 $y^k_j=x^k_j+\beta_kv^k_j>0$.
\STATE \hspace{0.6cm}If $\phi(y^k)\leq\phi(x^k)$ and inequality (\ref{ieq_1}) holds. \hspace{2cm}($\ast$)\\
\STATE
       \hspace{0.9cm}logic=false, $x^{k+1}=P(y^k)$, $\beta_{k+1}=\beta_k$, $k=k+1$.\\
\hspace{0.6cm}end(if)\\
\hspace{0.6cm}$\beta_k=\gamma\beta_k$. \\
 ~ end(while)
\end{algorithmic}
\end{algorithm}

Because the perturbation direction $v^k$ is selected as the
decreasing direction of $\phi$ at $x^k$,
 the size of $\beta_k$ represents the strength of regularization in some sense.
 Because the condition 3 of theorem \ref{ThEM} is very strict, the numerical experiments show that
  $\{\beta_k\}$  goes to zero very fast, which results in the regularization is very weak.
 Therefore,
we propose a modified version of algorithm \ref{algm1}, which is shown in algorithm \ref{algm2}. In algorithm \ref{algm2}, we only validate $I_A^b(x^{k+1})<I_A^b(x^k)$, rather than the inequality (\ref{ieq_1}) of theorem \ref{ThEM}. Since the inequality (\ref{ieq_1}) implies $I_A^b(x^{k+1})<I_A^b(x^k)$ by theorem \ref{ThEM},
algorithm 2 can be seen as a relaxation  of algorithm \ref{algm1}.
Furthermore, we introduce a relative decrease of $I_A^b(\cdot)$ to avoid the
situation that the amount of $I_A^b(x^k)-I_A^b(x^{k+1})$ for each iteration
is too small, which can accelerate the convergence of $I_A^b(x^k)$ and $x^k$ intuitively.

\begin{algorithm}[h]
\caption{Modification of the $\phi$-superiorization algorithm \ref{algm1}} %算法的标题
\label{algm2} %给算法一个标签，这样方便在文中对算法的引用
\begin{algorithmic}
\STATE $\beta_0>0$, $x^0>0$, $k=0$, and $0<\gamma<1$.
\STATE repeat: logic=true
\STATE \indent~~while logic
\STATE \indent~~~~ find a decreasing direction $v^k$ of $\phi$ at $x^k$, such that
 $y^k_j=v^k_j+\beta_kv^k_j>0$.
\STATE \hspace{0.6cm}If $\phi(y^k)\leq\phi(x^k)$ and $I_A^b(Py^k)<I_A^b(x^k)$.\hspace{2cm}($\star$)\\
\STATE
       \hspace{0.9cm}logic=false, $x^{k+1}=P(y^k)$.\\% $\beta_{k+1}=\beta_k$,\\
\STATE
  \hspace{0.9cm}If $\frac{I_A^b(x^k)-I_A^b(x^{k+1})}{I_A^b(x^k)}<Q_1$.\\%\hspace{2cm}($a_1$)\\
     \hspace{1.1cm} $\beta_{k+1}=\gamma\beta_k.$                 \\%\hspace{3.2cm} ($a_2$)\\
  \hspace{0.9cm}else\\
     \hspace{1.1cm} $\beta_{k+1}=\beta_k$.\\
  \hspace{0.9cm}end(else)\\
  \hspace{0.9cm} $k=k+1$
\STATE \hspace{0.6cm}else\\
\hspace{0.9cm}
 $\beta_k=\gamma\beta_k$                                      \\%\hspace{3.9cm}($a_3$)\\
 \indent \hspace{0.6cm}end(else)\\
 ~ end(while)
\end{algorithmic}
\end{algorithm}

\indent
In order to confirm the inequality (\ref{ieq_1}) of theorem \ref{ThEM} in algorithm \ref{algm1}, we should compute
$\sum_{j\in{S_k^-}}\hat{x}_j$ and $\sum_{j\in{S_k^+}}\hat{x}_j$ for each iteration. However,  $\hat{x}$ is unknown in the iterative procedure. In practice, we estimate the
values $\sum_{j\in{S_k^-}}\hat{x}_j$ and $\sum_{j\in{S_k^+}}\hat{x}_j$ by
$
\sum_{j\in{S_k^-}}\hat{x}_j\approx\frac{|S_k^-|}{N}B,
$
and
$
\sum_{j\in{S_k^+}}\hat{x}_j\approx\frac{|S_k^+|}{N}B, \label{est2}
$,
where $|S_k^-|$ and $|S_k^+|$ denote the cardinalities of the sets $S_k^-$ and $S_k^+$, respectively. Here, the capital letters $N$ and $B=\sum_{i=1}^Mb_i$ denote the length of $x$ and the total counts. Lastly, the parameter
$\rho$ is  chosen as $10^{-6}$ in this paper.
\\
\indent
In the following, we will discuss how to generate desirable perturbation $\beta_kv^k$ or $y^k$ for two concrete objective functions, TV and $l^1$-norm,
 such that the perturbations satisfy the conditions ($\ast$) and ($\star$) of
the two $\phi$-superiorization algorithms.
Since the TV regularization allows the reconstructed image to have sharp edges,
TV based models are widely used in imaging sciences \cite{Pan2006,pan2008,Censor2009,ROF1992}.
For an $H\times W$ image $x$ whose pixel values are denoted by $x_{i,j}$, the TV of $x$ is defined as
\begin{eqnarray}
TV(x)&=&\sum_{i=1}^{H-1}\sum_{j=1}^{W-1}\sqrt{(x_{i+1,j}-x_{i,j})^2+(x_{i,j+1}-x_{i,j})^2},\label{EqTv1}
\end{eqnarray}
where $H,W$ are the height and width of $x$. In order to reduce the value of TV at $x^k$, we choose $v^k$ as
\begin{eqnarray}
v^k=s^k/|s^k|_\infty
\end{eqnarray}
where $s^k\in\partial TV(x^k)$ is the sub-gradient of $TV$ at point $x^k$, and $|s^k|_\infty$ is the
maximum absolute value of the components of $s^k$. Therefore, the sequence $\{v^k\}$ is bounded. In fact, $v^k$ is the normalization of $s^k$
 in the $l^\infty$ space, rather than in the $l^2$ space used in \cite{Censor_2010,Herman2012}. \\
\indent
 In addition to the TV based models,
$l^1$-norm minimization method is another widely used technique in image sciences
\cite{Shen2012,Herman2011,Candes2006,Wang2009}.
Here, the $l^1$-norm is about the wavelet coefficients of $x$, which
is defined as
\begin{eqnarray}
\|T_{\{\psi_j\}}(x)\|_1&=&\sum_{j=1}^N|\alpha_j|,\label{EqL1}%\quad\textrm{with}\quad x=\sum_k\beta_k\psi_k(x)
\end{eqnarray}
where $\alpha_j$s  are the coefficients of $x$ under a given wavelet basis $\{\psi_j\}$,
and the letter $T$ denotes the wavelet decomposition operator. Although we can use the same
method for TV function to reduce the wavelet $l^1$-norm,
we introduce two more effective methods, soft and hard thresholding schemes,
to reduce the wavelet $l^1$-norm of $x^k$. Let
\begin{eqnarray}
&\gamma_j^k=&\text{Hard}(\alpha^k_j)=\left\{\begin{array}{ll}
0          &|\alpha_j^k|\geq{\beta_k}\\
-\frac{\alpha_j^k}{\beta_k}\text{sign}(\alpha_j^k)&|\alpha_j^k|<\beta_k\end{array}\right.,\\
\text{or}\notag\\
&\gamma_j^k=&\text{Soft}(\alpha_j^k)=\left\{\begin{array}{ll}
-\text{sign}(\alpha^k_j)&|\alpha_j^k|\geq{\beta_k}\\
-\frac{\alpha_j^k}{\beta_k}\text{sign}(\alpha^k_j)                         &|\alpha_j^k|<\beta_k\end{array}\right.,
\end{eqnarray}
where $\alpha_j^k$ are the wavelet coefficients of $x^k$. The perturbation direction for each iteration is defined as $v^k=T^{-1}_{\{\psi_j\}}(\gamma_j^k)$, and $y^k=x^k+\beta_kv^k$.
In fact, we need not to compute $v^k$ explicitly. By the linearity of wavelet transform, we can obtain the wavelet coefficients of $y^k$ directly by
\begin{eqnarray}
&\alpha_j^k=&\text{Hard}^C(\alpha^k_j)=\left\{\begin{array}{ll}\alpha_j^k&|\alpha_j^k|\geq{\beta_k}\\0&|\alpha_j^k|<\beta_k\end{array}\right.,\\
\text{or}\notag\\
&\alpha_j^k=&\text{Soft}^C(\alpha_j^k)=\left\{\begin{array}{ll}\alpha_j^k-\text{sign}(\alpha^k_j)\cdot\beta_k&|\alpha_j^k|\geq{\beta_k}\\0&|\alpha_j^k|<\beta_k\end{array}\right.,
\end{eqnarray}
and $y^k=T_{\{\psi_j\}}^{-1}(\alpha_j^k)$. In the numerical experiments,
we use Daubechies 6.8 bi-orthogonal wavelets with symmetric extensions at the boundaries. For referred convenience, we use hard-superiorization and soft-superiorization to distinguish them
for $l^1$-superiorization of EM algorithm in this paper.
%%%%%%%%%%%%
%And $\gamma$ is chosen as $0.5$
Furthermore, in order to avoid $y^k_j\leq0$, $y_j^k$ is
set as ${1\over2}x_j^k$ if $y_j^k\leq0$ in practice.
\section{Numerical Results}\label{Num_Res}
\begin{figure}
\centering
\begin{tabular}{ccc}
\includegraphics[width=5cm]{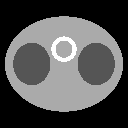}&
\includegraphics[width=5cm]{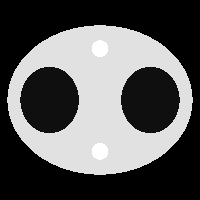}
\end{tabular}
\caption{activity map(left) and attenuation map(right). }\label{Figact_att}
\end{figure}

In this section we investigate the performance of the proposed algorithms
by several numerical experiments of SPECT reconstruction.
To this end,
and the projection data were generated based on the following model.
As shown in figure \ref{Figact_att}, the activity phantom
consists of an ellipsoidal background (body region) with axes of length 22.5cm and 30cm, which contains
two smaller ellipsoidal regions(lungs) with axes of length 10cm and 8.8cm, and  a ring(myocardium) of inner and outer diameters 6cm and  8cm, respectively.
The activities in myocardium, background, and lungs are
specified to be in the ratio 3:2:1.
\\
\indent
To simulate the attenuation coefficient in chest, we utilized the phantom
used in {\cite{Osem1994}}, which imitates a section of human thorax.
Besides the body background and lungs, the attenuation map consists of
two circular regions (bones) of diameter 2.5cm(see figure 1).
The attenuation coefficients were
0.03cm$^{-1}$ within 'lung' regions, 0.17cm$^{-1}$ within 'bone' regions, 0.15cm$^{-1}$ elsewhere
within the body ellipse, and 0.00cm$^{-1}$ outside the body.
\\
\indent
In the experiments, the activity and attenuation maps were evenly sampled in
 $[-15,15]\times[-15,15]$ on a grid of $128\times128$.
A perfect parallel hole collimator was assumed, and noise-free projection data were created via attenuated Radon transform formula (\ref{eq1}), which included tissue attenuation, but neither scattering nor blurring. In order to simulate the quantum noise in the simulated data, the following procedure was implemented \cite{fessler}.
\begin{itemize}
\item The projection data are scaled (multiplied by a constant factor)
so that the number of counts is a predefined integer.
\item Each value in the data set is then replaced by a random realization of a Poisson
variant with a mean equal to that value.
\end{itemize}
\indent
Two data sets, sixty and thirty projections, were generated over 180$^\circ$ evenly with view angles $\varphi_l=\frac{l-1}{N_0}\pi(l=1,\cdots,N_0$, and $N_0=60~ \text{or} ~30$). The counts were recorded in 128 bins per projection, and the
total counts were approximately 500K and 100K for two projection data sets, respectively.
\begin{figure}
\centering
\begin{tabular}{cccc}
\includegraphics[width=4cm]{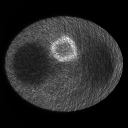}&
\includegraphics[width=4cm]{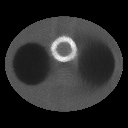}\\
\includegraphics[width=4cm]{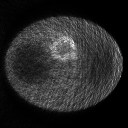}&
\includegraphics[width=4cm]{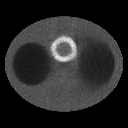}
\end{tabular}
\caption{Images reconstructed from the simulated data sets 1(top row) and 2(bottom row) by the classic EM algorithm(left column) and the mean images of 100 trials(right column).}\label{fig-em}
\end{figure}

\indent
We illustrated four numerical experiments for the proposed algorithms.
The first and second experiments are about the projection data set 1(60 projections, 500K) and projection
data set 2(30 projections, 100K) to validate the performance of the algorithms for different counts level.
The third experiment uses the randomly initial image to test the robustness of the proposed algorithms.
In the fourth experiment, we relax the condition of the TV-superiorized EM algorithms 1 and 2 further.
\\
\indent In order to evaluate the qualities of the reconstructed images,
we computed the mean images
$x^{\ast}$ from 100 noise trials as the standard images of the two data sets, respectively.
\begin{eqnarray}
x^{\ast}={1\over 100}\sum_{m=1}^{100}\hat{x}^m,
\end{eqnarray}
where $\hat{x}^m(m=1,2,\cdots,100)$ is the reconstructed image
of the $m$th trial by the classic EM algorithm, and the number of iterations for each trial is 30.\\
%%%%%%%%%%%%%%
\indent As shown in figure \ref{fig-em},
 the mean images for the two data sets are clear visually, while
 those reconstructed by the classic EM algorithm are dominated by noise, especially the image for data set 2 due to the low count level.
Thus, we use the mean square error (MSE) between the estimation $x$ and the mean
image $x^\ast$ from the same data set to measure the image qualities,
where the MSE is computed by
\begin{eqnarray}
\text{MSE}(x)=\frac{1}{N}\sum_{j=1}^N(x_j-x_j^{\ast})^2.
\end{eqnarray}
 The outputs of the different algorithms were taken as the best estimations in terms of MSE.
 In order to compare the qualities of images reconstructed from different data sets, we introduced the relative MSE(RMSE) defined as
 \begin{eqnarray}
\text{RMSE}(x)=\sqrt{\frac{\sum_{j=1}^N(x_j-x_j^{\ast})^2}{\sum^N_{j=1}(x^\ast_{j})^2}},
\end{eqnarray}
\\
\indent
In the numerical experiments, an uniformly initial image $x^0$ of value $c={1\over N}\sum_{i=1}^Mb_i$
was used for all experiments, unless there was a further explanation. The parameters $\beta_0$ were chosen as $c/2$ and $c/10$
for the TV- and $l^1$-superiorized EM algorithms, respectively.
The thresholding $Q_1$ and
the parameter $\gamma$ were chosen as $0.01$ and $1/2$, respectively.
In the following, we use TV-, hard- and soft-alg $n$($n=1,2$) as the abbreviations of the TV-, hard- and soft-superiorized EM algorithm $n$($n=1,2$) for simplicity.
\\
{\bf Experiment 1: sixty projections and 500K counts }
\begin{figure}[h]
\centering
\begin{tabular}{cccc}
&\includegraphics[width=4cm]{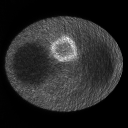}
&\includegraphics[width=4cm]{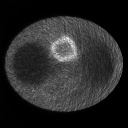}&
\includegraphics[width=4cm]{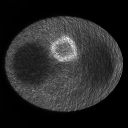}\\
&\includegraphics[width=4cm]{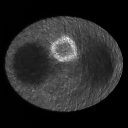}
&\includegraphics[width=4cm]{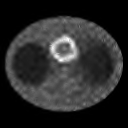}&
\includegraphics[width=4cm]{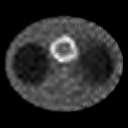}
\end{tabular}
\caption{
Reconstructed images of experiment 1. The images in top and bottom rows are reconstructed by the superiorized EM algorithms \ref{algm1} and \ref{algm2}, respectively.
The images reconstructed by TV-, hard- and soft-superiorized EM algorithms are displayed from column 1 to column 3.
}\label{fig-simul1}
\end{figure}
\begin{figure}[h]
\centering
\includegraphics[width=14cm]{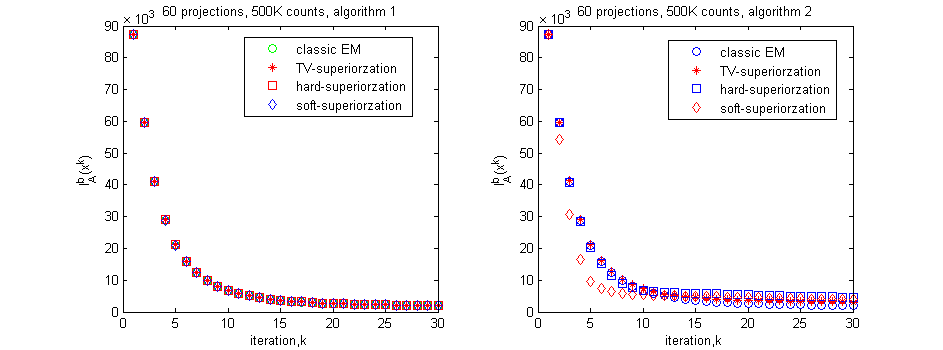}\\
\includegraphics[width=14cm]{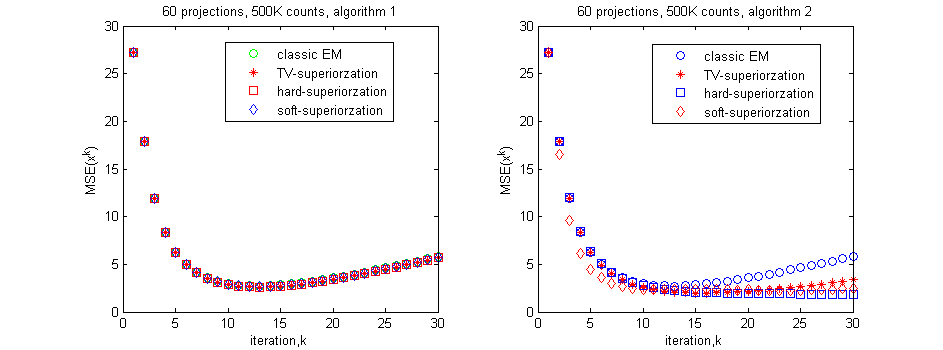}\\
\includegraphics[width=14cm]{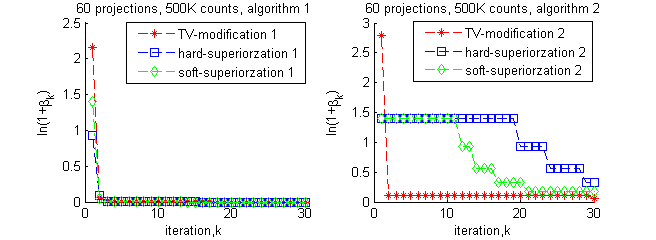}
\caption{Evolutions of $I_A^b(x^k)$(top row), MSE$(x^k)$(middle row) and $\ln(1+\beta_k)$(bottom row) of experiment 1.}\label{fig-simul11}
\end{figure}
\begin{table}[h]
\centering
\caption{TV values, $l^1$-norms, RMSEs and iterations of the reconstructed images of experiment 1.}\label{Table-simul1}
\begin{tabular}{ccccc}
\hline\hline
& EM& TV-alg 1&hard-alg 1&soft-alg 1\\
\hline
TV($\times10^3$)   & 26.708 &25.318   & 25.577  &  25.958\\%1.3987e+004
$l^1$($\times10^3$)&  12.953&12.404   & 12.320  &   12.489
\\%8.5244e+003
RMSE&0.1914&0.1873  &  0.1875  &  0.1888\\
iteration &13&13  &  13   & 13\\
\hline
& mean&TV-alg 2&hard-alg 2&soft-alg 2\\
\hline
TV($\times10^3$)    & 13.987&  19.434 &   14.115  &  15.443\\%7.1677e+003
$l^1$($\times10^3$) & 8.5244& 10.279   &  2.867   &  3.189\\%4.0374e+003
RMSE&-& 0.1633  &  0.1563   & 0.1747\\
iteration& -&15  &  30  &  23\\
\hline\hline
\end{tabular}
\end{table}

As expected, the images by the superiorized EM algorithms
are superior to the one by the classic EM algorithm in terms of TV value, $l^1$-norm and RMSE
(see figure \ref{fig-simul1} and table \ref{Table-simul1}), though the images by the superiorized algorithm 1 are
visually indistinguishable from the one by the classic EM algorithm.
Furthermore, the superiorized EM algorithm \ref{algm2} are superior to the
superiorized EM algorithm \ref{algm1},
because the regularization parameter $\beta_k$ goes to zero very fast for the superiorized EM algorithm \ref{algm1}.
\\
\indent As figure \ref{fig-simul11} shown, the evolutions of $I_A^b(x^k)$ and MSE$(x^k)$ of the
 superiorized EM algorithms 1 validate the conclusions of theorem \ref{ThEM}. Furthermore,
the evolutions of $I_A^b(x^k)$ and MSE$(x^k)$ of show the convergence of the
 superiorized EM algorithms 2, although we can not prove it theoretically.
Because the inequality of condition 3 is very strict, the parameters $\beta_k$ went to
zero very fast for the superiorized algorithm \ref{algm1}. This results in the reconstructed images and the evolutions of $I_A^b(x^k)$ and MSE$(x^k)$
of the superiorized EM algorithms \ref{algm1} and the classic EM algorithm are indistinguishable from each other,
i.e. the regularization of the superiorized EM algorithm 1 is not remarkable.
\\
\indent
%As shown in table \ref{Table-simul1},
Although the number of iterations
of the superiorized EM algorithm 2 is larger than the classic EM and the
superiorized EM algorithm 1(table \ref{Table-simul1}), the superiorized algorithm 2 is superior
to the superiorized algorithm 1 in terms of RMSE: RMSE$(x^{tv,13})=0.1668$,
RMSE$(x^{h,13})=0.1741$ and RMSE$(x^{s,13})=0.1797$,
where $x^{tv,13},x^{h,13},x^{s,13}$ denote the results of tv-, hard- and soft-superiorized algorithms 2, respectively.
\\
{\bf Experiment 2: thirty projections and 100K counts }
\\
\begin{figure}[h]
\centering
\begin{tabular}{cccc}
&\includegraphics[width=4cm]{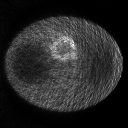}
&\includegraphics[width=4cm]{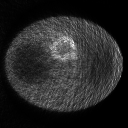}&
\includegraphics[width=4cm]{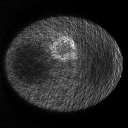}\\
&\includegraphics[width=4cm]{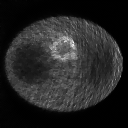}
&\includegraphics[width=4cm]{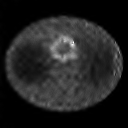}&
\includegraphics[width=4cm]{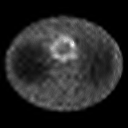}\\
\end{tabular}
\caption{Reconstructed images of experiment 2. The images reconstructed
by the TV-, hard- and soft-superiorized EM algorithms 1 and 2 are displayed in
columns 2, 3 and 4, and rows 1 and 2, respectively. }\label{fig-simul2}
\label{noised(30)}
\end{figure}
%%%%%%%%%%%%%%%%%%%%%%
\begin{table}[h]
\centering
\caption{TV, $l^1$-norm, RMSE and iterations of the images
of experiment 2.}\label{Table-simul2}
\begin{tabular}{ccccc}
\hline\hline
  &EM&TV-alg 1&hard-alg 1&soft-alg 1\\
\hline
TV($\times10^3$)     & 13.739&  13.105     &  13.396  &  13.061\\%1.3987e+004
$l^1$($\times10^3$)  & 6.729 &  6.452     &  6.451   &  6.277
\\%8.5244e+003
RMSE &0.2822& 0.2778  &  0.2802 &   0.2795\\
iteration & 8&8&8&8\\
\hline
   &mean&TV-alg 2&hard-alg 2&soft-alg 2\\
\hline
TV($\times10^3$)     &7.167 & 9.769  &  5.733 &   5.576\\%7.1677e+003
$l^1$($\times10^3$)  & 4.037 & 5.069  &  1.434 &    1.283\\%4.0374e+0034.7527    1.4136    1.2625
RMSE &0& 0.2490  &  0.2091  &  0.2064\\
Iteration & &9&12&10\\
\hline\hline
\end{tabular}
\end{table}
\begin{figure}[h]
\centering
\includegraphics[width=14cm]{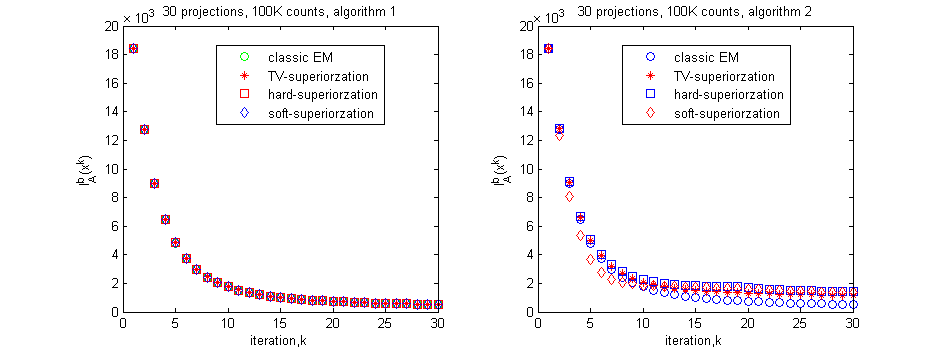}\\
\includegraphics[width=14cm]{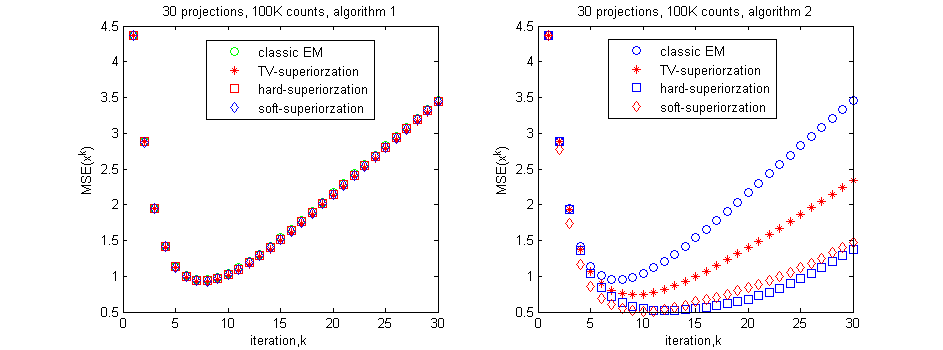}\\
\includegraphics[width=14cm]{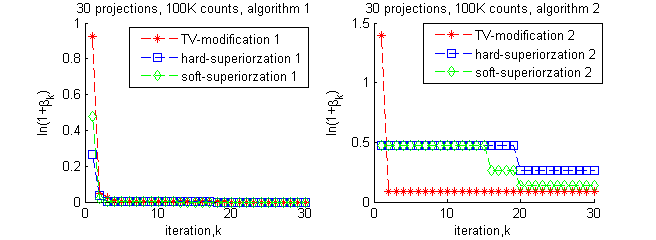}
\caption{Evolutions of experiment 2. The evolutions of $I_A^b(x^k)$(top row),  MSE$(x^k)$(middle row) and
$\ln(1+\beta_k)$(bottom row) of different algorithms.}\label{fig-simul21}
\end{figure}
\indent
From the observations of figures \ref{fig-simul2} and \ref{fig-simul21}, and table \ref{Table-simul2},
we can draw the same conclusions as the experiment 1.
Comparing figures \ref{fig-simul1} and \ref{fig-simul2},
we can see that the visual quality of the images reconstructed from data set 2
are inferior to those reconstructed from data set 1
because of the low count level.
And the performance of the superiorized EM algorithms \ref{algm2} are much more remarkable than the superiorized EM algorithms \ref{algm1} for data set 2.
\\
\indent
Comparing the results of the two experiments above, we can observe that
the TV-superiorized EM algorithm \ref{algm2} is better in terms of RMSE, while the $l^1$-superiorized EM algorithms \ref{algm2}
in terms of TV and $l^1$-norm values. In addition, the thresholding operations cause the Gibbs oscillations in the reconstructed images by $l^1$-superiorized EM algorithms \ref{algm2}.
\\
\indent
Because the detectors rotated from top to bottom through
left side, the gamma rays emitted from the right part pixels are more likely to
be absorbed. Therefore,
 the right part of the reconstructed images are blurred strongly(see figures \ref{fig-em},  \ref{fig-simul1} and \ref{fig-simul2}).
\\
%%%%%%%%%%%%%%%%%%%%%%%%%%%%%%%%%%%%%%%%%%%%%%%%%%%%%%%%%%%%
{\bf Experiment 3: initial image $x^0$ with random values on interval $[1,2]$ for data set 1}
\\
\indent
In this experiment, an initial image $x^0$ with random values on interval $[1,2]$ and the projection data set 1 are used.
The reconstructed images by different algorithms are displayed in figure \ref{fig-simul4},
and the corresponding TV values, $l^1$-norms and RMSEs are tabulated in table \ref{Table-simul4}.
Because the evolutions of $I_A^b(x^k)$ and MSE$(x^k)$ are very similar to these of experiment 1, we only plot the evolution of
$\ln(1+\beta_k)$ in figure \ref{fig-simul41}.
\\
\indent
This experiment shows that the superiorized EM algorithms 2 are stable and robust for initial image.
By comparing figures \ref{fig-simul1} and \ref{fig-simul4}, and tables \ref{Table-simul1} and \ref{Table-simul4}, we have that the effect of initial image is very strong to the
classic EM algorithm and the superiorized EM algorithm \ref{algm1},
but very weak to the superiorized EM algorithm \ref{algm2}.
\\
\indent
A surprising observation is that the randomly initial image is superior to
the uniformly initial image for the superiorized EM algorithms 2 in term of RMSE by comparing table \ref{Table-simul1} and table \ref{Table-simul4}. This changes the long-standing opinion about the
selection of the initial image for EM-like algorithm, and present a new method to improve the qualities
of the reconstructed image.
\begin{table}[h]
\centering
\caption{TV values, $l^1$-norms, RMSEs and iterations of the images
of experiment 3.}\label{Table-simul4}
\begin{tabular}{ccccc}
\hline\hline
  &EM&TV-alg 1&hard-alg 1&soft-alg 1\\
\hline
TV($\times10^3$)     & 40.591 &   29.804  &  40.364  &  34.214\\%1.3987e+004
$l^1$($\times10^3$)  & 19.098 &   14.632  &  18.907  &  16.138
\\
RMSE &0.2538&0.2112  &  0.2529 &   0.2266\\
iteration & 13 & 13 &13   & 13\\
\hline
  & -&TV-alg 2&hard-alg 2&soft-alg 2\\
\hline
TV($\times10^3$)     &-&17.133 &   16.203   & 16.372\\%7.1677e+003
$l^1$($\times10^3$)  & -  &  9.403 &   3.417  &  3.207\\%4.0374e+0034.7527    1.4136    1.2625
RMSE &-&0.1469   & 0.1683  &  0.1873\\
Iteration &-& 19  &  30 &   23\\
\hline\hline
\end{tabular}
\end{table}
\begin{figure}[h]
\centering
\begin{tabular}{cccc}
&\includegraphics[width=4cm]{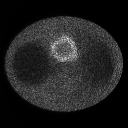}\\
\includegraphics[width=4cm]{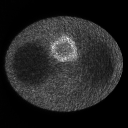}
&\includegraphics[width=4cm]{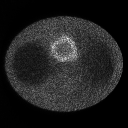}
&\includegraphics[width=4cm]{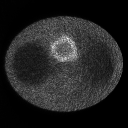}\\
\includegraphics[width=4cm]{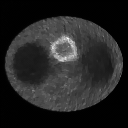}
&\includegraphics[width=4cm]{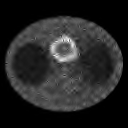}&
\includegraphics[width=4cm]{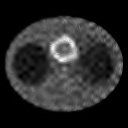}
\end{tabular}
\caption{Images of experiment 3. The top image is reconstructed by the classic EM algorithm.
The images in columns 1, 2 and 3 and the middle and bottom rows  are reconstructed
by the TV-, hard- and soft-superiorized EM algorithms 1 and 2, respectively. }\label{fig-simul4}
\end{figure}
\begin{figure}[h]
\centering
\includegraphics[width=14cm]{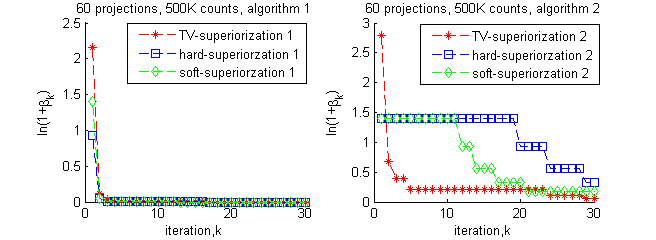}
\caption{ Evolutions of $\ln(1+\beta_k)$ of experiment 3.}\label{fig-simul41}
\end{figure}

\begin{figure}
\centering
\begin{tabular}{cc}
\includegraphics[width=4cm]{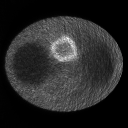}
% tv :25.736 10^3, l^1: 12.320
&\includegraphics[width=4cm]{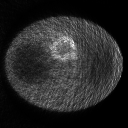}\\
\includegraphics[width=4cm]{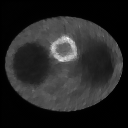}
%tv: 14.432 10^3,l^1: 10.765
&\includegraphics[width=4cm]{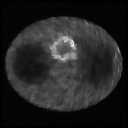}
\end{tabular}
\caption{Images of experiment 4. The images in left and right columns are reconstructed from data sets 1 and 2, while the images in first and second rows are reconstructed by modified versions of TV-superiorized EM algorithms 1 and 2, respectively. From left to right and from top to bottom, the images are labeled as image 1 to image 4 for reference.}\label{fig-simul3}
\end{figure}

\begin{figure}[h]
\centering
\includegraphics[width=14cm]{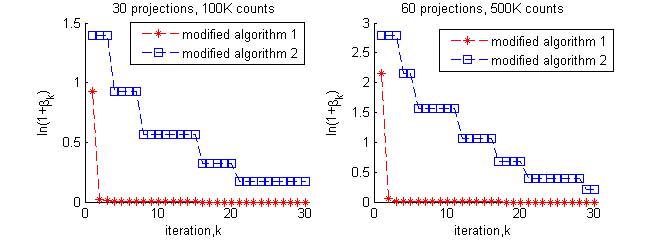}
\caption{Evolutions of $\ln(1+\beta_k)$  of experiment 4. Left for
data set 1 and right for data set 2, respectively.}\label{fig-simul31}
\end{figure}

\indent
Comparing the evolutions of $\ln(1+\beta_k)$ in figures \ref{fig-simul11}, \ref{fig-simul21} and \ref{fig-simul41}, we can observe the following facts. Firstly, the parameter $\beta_k$ for the superiorized EM algorithms \ref{algm1} go to zero fast because of the strict condition 3 of theorem \ref{ThEM}.
Secondly, the
reconstructed images of uniformly initial image is more smooth
than these of randomly initial image by the TV-superiorized algorithm at low iterations,
which results in the parameter $\beta_k$ decrease at fast and low rates for experiments 1 and 3,
respectively,  because of the decreasing condition of the TV function in the conditions ($\ast$) and ($\star$).
Thirdly, the thresholding operations(hard and soft) always reduce the $l^1$-norm, which implies that the
conditions ($\ast$) and ($\star$) for the $l^1$-superiorized algorithm 1 and  2 become  one condition about the decreasing of K-L distance.
\\
\indent The last observation enlightens us to modify the TV-superiorized EM algorithm further, which  discards
the deceasing condition of TV function in the conditions ($\ast$) and ($\star$).\\
%%%%%%%%%%%%%%%%%%%%%%%%%%%%%%%%%%%%%%%%%%%%%%%%%%%%%%%%%
{\bf Experiment 4: modification of the TV-superiorized algorithms 1 and 2}
\\
\indent
Figure \ref{fig-simul3} displays the reconstructed images by the modified versions of TV-superiorized EM algorithms 1 and 2 in absence of the
decreasing condition of TV function.
The TV values, $l^1$-norms, RMSEs and iterations are tabulated in table \ref{Table-simul3}. And
figure \ref{fig-simul31} plots the evolutions of $\ln(1+\beta_k)$.
\\
\indent As our expectation, the evolution of parameter $\beta_k$ of the modified version of TV-superiorized algorithm 2 is similar to the  $l^1$-superiorized EM algorithm 2, decreasing at much slower rate.
However, this modification has very little effect on the TV-superiorized algorithm 1.
\\
\indent
It is amazing that the modified algorithms also reduce the TV function (see table \ref{Table-simul4}), even the reconstructed images are better than those reconstructed by  the TV-superiorized EM algorithm \ref{algm2},
although we do not validate the decreasing condition of it. The reasons include two aspects.
In superiorized algorithms 1 and 2, there are two conditions to control the decreasing of $\beta_k$, which causes the size of $\beta_k$  is very small at large iterations. Therefore, the strength of regularization is very weak, and the reconstructed image is not good enough. For the modified superiorization algorithms, there is only one condition to control the decreasing of $\beta_k$, and $\beta_k$ decreases at a lower rate. Therefore, the modified superiorization algorithms can maintain stronger regularization at large iterations, and the reconstructed image is much better.
The further study about this algorithm is future work.

\begin{table}[h]
\centering
\caption{TV values, $l^1$-norm, RMSEs and iterations of the images of experiment 4. }\label{Table-simul3}
\begin{tabular}{ccccc}
\hline\hline
  & image 1& image 2& image 3& image 4\\
\hline
%%30:
TV($\times10^3$) & 25.318 & 12.931  & 12.955 &   7.080\\%13.105  &  7.758  \\
$l^1$($\times10^3$) &  12.404  &  6.408 & 7.013  &  3.399  \\
%12.534  &  7.494  & 6.4517  &  3.6920 \\
RMSE &0.1873   & 0.2772  &0.1231   &0.1999 \\
iteration &13 & 8  &  29   & 18\\
\hline\hline
\end{tabular}
\end{table}
\begin{Remark}
The experiments show that the superiorized EM algorithm 2 is convergent, though we cannot prove it theoretically so far.
\end{Remark}
\begin{Remark}
The parameter $\beta_k$ represents the strength of regularization in a sense,
so we can obtain the regularization reconstruction by terminating the algorithms
as long as $\beta_k$ is smaller than a predefined threshold.
 As explanation above, the size of $\beta_k$ represents the strength of regularization. Therefore, the threshold  should be related to the noise level.
 Intuitively, in order to maintain the regularization strength, the higher the noise level is, the larger the threshold is.
The further discussion about the selection of it will be studied in future work.
\end{Remark}
\section{Conclusions and Discussions}\label{discussion}
In this paper, the convergence of the EM algorithm in the presence
of perturbations is discussed, and
the superiorized EM algorithm
based on the convergent conditions and its modified version
are proposed.
The numerical experiments validate the correction of theorem \ref{ThEM}.
The superiorized EM algorithms could
efficiently reduce the corresponding objective functions which we are superiorizing. Furthermore,
The proposed algorithms are more stable and robust than the classic EM algorithm for low counts projection data and randomly initial image.
\\
\indent
Although the numerical experiments show the convergence of algorithm 2,
 we cannot prove the convergence of it theoretically.
A more challenging work is about the amazing observation of the experiment 4, which enlightens us
to modify the superiorized EM algorithm further.
In addition, we could not prove $\phi(x^\ast)\leq \phi(\hat{x})$ theoretically,
where $\hat{x}$ and $x^\ast$ are
the solutions by the classic iteration algorithm and the $\phi$-superiorization
version \cite{Censor_2010}.
\section*{Acknowledgement}\label{appendix}
This work is supported by the
National Basic Research Program of China (2011CB809105) and NSF grants of China (61121002, 10990013).
The authors are grateful for the helpful discussions with Professor Haomin Zhou(School of Mathematics,
Georgia Institute of Technology).
\bibliographystyle{ieeetr}
\bibliography{spect_superiorization49}
\end{document}